\documentclass[a4paper,UKenglish,cleveref, autoref, thm-restate]{lipics-v2021}
\nolinenumbers

\usepackage{hyperref}

\usepackage{amsmath}
\usepackage{upgreek}
\usepackage{xcolor}
\usepackage{colortbl}
\usepackage[capitalise]{cleveref}

\usepackage{pbox}
\usepackage{makecell}

\usepackage{tikz-cd}
\usepackage{tabularx}
\definecolor{keywordcolor}{rgb}{0.7, 0.1, 0.1}   
\definecolor{tacticcolor}{rgb}{0.0, 0.1, 0.3}    
\definecolor{commentcolor}{rgb}{0.4, 0.4, 0.4}   
\definecolor{stringcolor}{rgb}{0.5, 0.3, 0.2}    
\definecolor{symbolcolor}{rgb}{0.1, 0.2, 0.7}    
\definecolor{sortcolor}{rgb}{0.1, 0.5, 0.1}      
\definecolor{attributecolor}{rgb}{0.7, 0.1, 0.1} 
\definecolor{errorcolor}{rgb}{1, 0, 0}           

\usepackage{listings}
\usepackage{flushend}
\usepackage{xspace}
\usepackage{underscore}
\usepackage{microtype}
\usepackage{datetime}
\usepackage{MnSymbol} 
\usepackage{ragged2e}
\newtheorem{assumption}{Assumption}
\newcommand{\lean}[1]{\lstinline[language=lean]{#1}}

\usepackage{scalefnt}
\title{Formalization of Optimality Conditions for Smooth Constrained Optimization Problems} 

\author{Chenyi Li}{School of Mathematical Sciences, Peking University, China}{lichenyi@stu.pku.edu.cn}{}{} 
\author{Shengyang Xu}{School of Mathematical Sciences, Peking University, China}{xushengyang0429@gmail.com}{}{}
\author{Chumin Sun}{Huawei Technologies Co., Ltd.}{sunchumin@huawei.com}{}{}
\author{Li Zhou}{Huawei Technologies Co., Ltd.}{zhouli107@huawei.com}{}{}
\author{Zaiwen Wen\footnote{corresponding author}}{Beijing International Center for Mathematical Research, Peking University, China}{wenzw@pku.edu.cn}{}{}

\authorrunning{C.\,Li, S.\,Xu, C.\,Sun, L.\,Zhou, and Z.\,Wen}

\Copyright{Chenyi Li, Shengyang Xu, Chumin Sun, Li Zhou, and Zaiwen Wen} 



\category{} 

\relatedversion{} 

\supplement{Our code for this paper can be found at \url{https://github.com/optsuite/optlib}}

\acknowledgements{We want to thank the students from the summer school and winter school at Peking University.}

\begin{document}

\maketitle

\begin{abstract}
Optimality conditions are central to analysis of optimization problems, characterizing necessary criteria for local minima. Formalizing the optimality conditions within the type-theory-based proof assistant Lean4 provides a precise, robust, and reusable framework essential for rigorous verification in optimization theory. In this paper, we introduce a formalization of the first-order optimality conditions (also known as the Karush-Kuhn-Tucker (KKT) conditions) for smooth constrained optimization problems by beginning with concepts such as the Lagrangian function and constraint qualifications. The geometric optimality conditions are then formalized, offering insights into local minima through tangent cones. We also establish the critical equivalence between the tangent cone and linearized feasible directions under appropriate constraint qualifications. Building on these key elements, the formalization concludes the KKT conditions through the proof of the Farkas lemma. Additionally, this study provides a formalization of the dual problem and the weak duality property.
\end{abstract}

\section{Introduction}
The first-order optimality conditions (also known as the Karush-Kuhn-Tucker (KKT) conditions) are fundamental in optimization theory, providing a framework for analyzing and understanding the behavior of solutions in smooth constrained optimization problems. These conditions provide a necessary criterion for optimality in smooth constrained problems, given appropriate constraint qualifications. The KKT conditions are described by Kuhn and Tucker in \cite{KuhnTucker+1951+481+492}. Mangasarian provides a thorough description of the constraint qualifications in \cite{Mangasarian1994Nonlinear}. Lagrange multipliers and optimality conditions for general problems are explored by Rockafellar in \cite{Rockafellar1993Lagrange}. Several analyses have examined the KKT conditions in various contexts. Achtziger et al. investigate mathematical programs with vanishing constraints \cite{achtziger2008mathematical}. The KKT conditions and constraint qualifications have been explored in bilevel optimization \cite{ye2006constraint}, semi-definite programming \cite{ye2008first}, and sequential quadratic programming \cite{qi2000constant}. The intrinsic formulation of KKT conditions on smooth manifolds is studied in \cite{bergmann2019intrinsic}, while the analysis is extended to Banach spaces in \cite{borgens2020new}. 

Numerical optimization problems have been studied in formalization in several literature, using different formalization languages, such as Coq \cite{huet1997coq}, Isabelle \cite{Nipkow2002APA}, and Lean \cite{de2015lean}. Some current research utilizes Lean4 to formalize optimization problems. CVXLean \cite{Bentkamp2023Verified} is introduced to perform disciplined convex programming on optimization problems through Lean. It verifies the reduction in a formalized way. Several manual and error-prone steps through formalization are discussed in \cite{Fernandez2024Transforming} to transform optimization problems into an equivalent form that is accepted by disciplined convex programming frameworks. To the best of our knowledge, there has been limited research on the formalization of optimality conditions for constrained optimization problems. Some studies formalize basic results in functional analysis and optimization. For example, the formalization of functional analysis with semilinear maps is explored by Dupuis et al. in \cite{dupuis_et_al:LIPIcs.ITP.2022.10}. Li et al. give the formalization of convergence rates for various first-order optimization methods \cite{li2024formalization}.  Doorn and Macbeth formalize the Gagliardo-Nirenberg-Sobolev inequality in \cite{vandoorn_et_al:LIPIcs.ITP.2024.37}. Moreover, the Doob’s martingale convergence theorems are formalized in Mathlib by Ying et al. \cite{Ying2023Doob}. 

In this paper, we present the formalization of the optimality conditions for smooth constrained optimization problems in Lean4 following Chapter 12 of \cite{nocedal1999numerical}. Our main contributions are as follows. Using the \lean{structure} type, we formally represent constrained optimization problems and introduce the foundational definitions, including feasible points and the Lagrangian function. To derive the optimality conditions, geometric optimality conditions are established based on the formalization of the tangent cone and the linearized feasible directions. The constraint qualifications, including the linearly independent constraint qualification (LICQ) and the linear constraint qualification, are formalized to prove the equivalence between the positive tangent cone and the linearized feasible cone at a given point. The Farkas lemma is formalized to leverage the information from the linearized feasible cone. We present a formalization for the proof of the KKT conditions based on these lemmas. Furthermore, we explore the dual problem and establish the weak duality property, illustrating the relationship between primal and dual formulations.

The remainder of this paper is organized as follows. In Section \ref{sec: lean-pre}, we provide an overview of Lean4 language and discuss existing formalizations of the positive tangent cone, which is frequently used throughout the paper. Discussion on the formalization of the constrained optimization problems is given in \ref{sec: problem}. The LICQ and the linear constraint qualification are presented in Section \ref{sec: CQ}. The formalization of the Farkas lemma is addressed in Section \ref{sec: Farkas}. The paper concludes with the formalization of the KKT conditions in Section \ref{sec: KKT}. We provide formalization related to the dual problem and the duality property in Section \ref{sec: weak duality}.

\section{Lean preliminaries}\label{sec: lean-pre}
\subsection{A Brief Introduction to Lean}
Developed by Microsoft Research since 2013, Lean is a formal language based on dependent type theory, designed for mathematicians and computer scientists to write precise mathematical proofs and perform logical reasoning. In Lean, logical propositions are treated as types, and proofs are constructions of these types. The dependent type system, along with the ability to define inductive structures, allows the type checker to verify the correctness of propositions, ensuring the reliability of mathematical proofs. Its intuitive syntax enables the expression of complex concepts and proofs in a way that closely resembles natural mathematical language, reducing the barrier to learn formalization.

Mathlib4, the core mathematical library \cite{mathlibcommunity} in the Lean system, aims to establish a standardized and modular mathematical knowledge base, facilitating efficient verification and sharing. Independent of Lean's core structure, Mathlib is a comprehensive repository built on Lean’s foundational language, much like the relationship between Python and packages such as numpy and torch. With over 170,000 theorems and 80,000 definitions, Mathlib spans a wide range of topics, from basic to applied mathematics, significantly enhancing the verifiability and shareability of mathematical knowledge.

\subsection{Tangent Cone}
The definition of the tangent cone is crucial in the optimality condition of a constrained optimization problem. This is already formalized in Mathlib4. In this subsection, we review this definition and explain some basic grammars in Lean.

The tangent cone of a point $x$ to the set $\Omega$ is defined as below.
\begin{definition}\label{def: tangent cone}
    The vector \( d \) is said to be a tangent vector to \( \Omega \) at a point \( x \) if there exists a feasible sequence \( \{ z_k \} \) converging to \( x \) and a sequence of positive scalars \( \{ t_k \} \) with \( t_k \to 0 \) such that
    \[\lim_{k\to\infty}\frac{z_k-x}{t_k} = d.\]
    The set of all tangent vectors to $\Omega$ at $x$ is called the tangent cone and is denoted by $T_\Omega\{x\}$.
\end{definition}

The formalization of the tangent cone differs from the definition in \cite{nocedal1999numerical} by utilizing a sequence that tends to infinity. However, this approach is mathematically equivalent to the definition provided in Definition \ref{def: tangent cone}. In Lean, the tangent cone at $x$ is represented as below
\begin{lstlisting}
def posTangentConeAt {E : Type u} [NormedAddCommGroup E] [NormedSpace ℝ E]  
    (s : Set E) (x : E) : Set E :=
  { y : E | ∃ (c : ℕ → ℝ) (d : ℕ → E), (∀ᶠ n in atTop, x + d n ∈ s) ∧
    Tendsto c atTop atTop ∧ Tendsto (fun n => c n • d n) atTop (nhds y) }
\end{lstlisting}
The type of the definition is \lean{Set E}, with its value behind the term \lean{:=}. In this definition, different types of brackets are used to express distinct sets of assumptions. Curly brackets enclose assumptions that can be inferred from the others. Square brackets denote the instances of the objective variables, which are automatically verified when specific variables are applied to \lean{E}. Round brackets indicate the assumptions that need to be explicitly provided when reusing the definition.

The limit behavior is expressed using the concept of ``filter”, which can be understood as a collection of sets. Mapping relations between two filters are used to describe limitation process. The notation \lean{atTop} represents filter for the limit at infinity. Thus, the expression \lean{Tendsto c atTop atTop} indicates that $\lim_{n \to \infty} c_n = +\infty$. Additionally, the notation $\forall^f$ indicates that a property holds when $n$ becomes sufficiently large.

\section{Formalization of Constrained Optimization Problem}\label{sec: problem}
In this section, we consider the formalization of constrained optimization problems. A general formulation of these problems is
\begin{align}
    \min_{x\in\mathcal{X}}\ f(x), \quad\ \text{subject to}\ \begin{cases}
        c_i(x) = 0,&\quad i \in\mathcal{E},\\
        c_i(x) \geqslant 0,&\quad i \in\mathcal{I},\\
    \end{cases}
    \label{eq: constrained problem}
\end{align}
where the objective function $f(x) : \mathbb{R}^n \to \mathbb{R}$ and the constraints $c_i(x)  : \mathbb{R}^n \to \mathbb{R} $ are real-valued functions defined on $\mathbb{R}^n$. The domain of this problem is $\Omega$ as a subset of $\mathbb{R}^n$, i.e. $\mathcal{X} \subseteq \mathbb{R}^n$. $\mathcal{E}$ and $\mathcal{I}$ are two finite sets of indices. $c_i(x),\ i\in\mathcal{E}$ denote the equality constraints and $c_i(x),\ i\in\mathcal{I}$ are the inequality constraints. We formalize the framework with the type \lean{structure} as follows:
\begin{lstlisting}
structure Constrained_OptimizationProblem (E : Type _) (τ σ : Finset ℕ) :=
    (domain : Set E) (objective : E → ℝ)
    (equality_constraints : (i : ℕ) → E → ℝ)
    (inequality_constraints : (j : ℕ) → E → ℝ)
    (eq_ine_not_intersect : τ ∩ σ = ∅)
\end{lstlisting}
The term \lean{E} denotes the general objective space that we consider for our optimization problem. In formulation \eqref{eq: constrained problem}, it is the Euclidean space $\mathbb{R}^n$. However, for general cases, we may consider a general Hilbert space as our objective space. Moreover, $\tau$ and $\sigma$ correspond to the indices of equality constraints $\mathcal{E}$ and inequality constraints $\mathcal{I}$, respectively. It is worth noting that in our formalization we give a concrete type \lean{Finset ℕ} rather than \lean{Fintype _} for $\tau$ and $\sigma$. The reason is that when considering constraints, we often label them with natural numbers such as $c_0, c_1 \in \mathbb{N}$. Due to the equivalence of different definitions, we adopt the indices as subsets of natural numbers. Meanwhile, since the equality and inequality constraints are defined on different sets, there is actually no need for $\tau$ and $\sigma$ to be disjoint. We add the condition \lean{eq_ine_not_intersect} for the clarity of the exposition, especially for the definition of active set, which will be mentioned later. It is obvious that every constrained optimization problem can be formalized as above. In subsequent contents, we always use \lean{p} to represent \lean{Constrained\_OptimizationProblem E τ σ} and give the following definition.
\begin{lstlisting}
variable {E : Type *} {τ σ : Finset ℕ} 
variable {p : Constrained_OptimizationProblem E τ σ}
\end{lstlisting}
For a constrained optimization problem, we define the feasible set $\Omega$ to be the set of points $x$ that satisfy the constraints as
\begin{align*}
    \Omega = \{x\left.\right| x \in \mathcal{X},  c_i(x) = 0,\ i\in \mathcal{E};\ c_i(x) \geqslant 0,\ i\in \mathcal{I}\}.
\end{align*}

This is formalized in Lean as
\begin{lstlisting}
def FeasPoint (point : E) : Prop := point ∈ p.domain ∧
    (∀ i ∈ τ, p.equality_constraints i point = 0) ∧
    (∀ j ∈ σ, p.inequality_constraints j point ≥ 0)

def FeasSet : Set E := {point | p.FeasPoint point}
\end{lstlisting}

In constrained optimization problems, one of our primary concerns is the characterization of local minima. We formalize this using the term \lean{IsLocalMinOn} in Mathlib.
\begin{lstlisting}
def Local_Minima (point : E) : Prop :=
  (p.FeasPoint point) ∧ IsLocalMinOn p.objective p.FeasSet point
\end{lstlisting}
Since \lean{IsLocalMinOn f s a} means that $f(a) \leq f(x)$ for all $x\in s$ in some neighborhood of $a$, there is no constraint that $a \in s$. Hence, we put an extra constraint \lean{p.FeasPoint point} to enforce the feasibility.

Corresponding to problem \eqref{eq: constrained problem}, we can form the Lagrangian function as
\begin{align}
    \mathcal{L}(x,\lambda,\mu) = f(x) - \sum_{i \in \mathcal{E}} \lambda_i c_i(x) - \sum_{i \in \mathcal{I}} \mu_i c_i(x),
    \label{eq: lagrangian function}
\end{align}
where the multipliers $\lambda_i$ are for equality constraints and $\mu_i$ are for inequality constraints. There are some equivalent definitions of the Lagrangian function and we take the one with minus operation. As a preliminary to stating the necessary conditions, the formalization of the Lagrangian function \eqref{eq: lagrangian function} for a constrained optimization problem is straightforward, as
\begin{lstlisting}
def Lagrange_function := 
    fun (x : E) (lambda : τ → ℝ) (mu : σ → ℝ) ↦ p.objective x -
    Finset.sum univ fun i ↦ (lambda i) * p.equality_constraints i x - Finset.sum univ fun j ↦ (mu j) * p.inequality_constraints j x
\end{lstlisting}
The type $\tau \to \mathbb{R}$ denotes a function mapping from the type $\tau$ to $\mathbb{R}$, where $\tau$ denotes the subtype of $\mathbb{N}$, which can be expressed as \{\lean{i // i ∈ τ}\}. 

From the perspective of the tangent cone, it is easy to see that if $x^\ast$ is a local minimum of the problem, there is no feasible descending direction in $T_\Omega(x^\ast)$. This is exactly the geometric optimality condition, which gives the necessary conditions that must be satisfied by any solution points. This theorem is mathematically stated as follows.
\begin{theorem}\label{thm: geometric optimality}
    Let $x^\ast$ be a local minimizer of the constrained problem \eqref{eq: constrained problem}. If the objective function $f$ is differentiable at $x^\ast$, then
\begin{align}
d^\top \nabla f(x^\ast) \geqslant 0,\quad \forall d \in T_\Omega(x^\ast),
\label{eq: geometric 1}
\end{align}
which is equivalent to
\begin{align}
T_\Omega(x^\ast) \cap \{d\ |\ d^\top \nabla f(x^\ast) < 0\} = \varnothing.
\label{eq: geometric 2}
\end{align}
\end{theorem}

We formalize \eqref{eq: geometric 1} and \eqref{eq: geometric 2} separately in Lean as follows
\begin{lstlisting}
variable [NormedAddCommGroup E] [InnerProductSpace ℝ E] [CompleteSpace E]
lemma posTangentCone_localmin_inner_pos {f : E → ℝ} {loc : E} (hl : IsLocalMinOn f p.FeasSet loc)(hf : DifferentiableAt ℝ f loc) :
    ∀ v ∈ posTangentConeAt p.FeasSet loc,
    ⟪gradient f loc, v⟫_ℝ ≥ (0 : ℝ)

theorem local_minima_TangentCone' (loc : E)
    (hl : p.Local_Minima loc) (hf : Differentiable ℝ p.objective) :
    posTangentConeAt p.FeasSet loc ∩
    {d | ⟪gradient p.objective loc, d⟫_ℝ < (0 : ℝ)} = ∅
\end{lstlisting}
We assume that the objective space is a Hilbert space, which is given as assumptions in the \lean{variable} environment for subsequent theorems in this paper. This constraint mainly comes from the utilization of gradient, which requires the domain space to be a Hilbert space. 

\section{Constraint Qualifications}\label{sec: CQ}
The definition of tangent cone relies only on the geometry of $\Omega$, and it is not easy to provide a specific formula of it everywhere. Therefore, we introduce the definition of the active set and linearized feasible directions. They depend on the value of the constraint functions $c_i(x)$ and are easy to characterize at a certain point. First, we need to give the definition of the active set.
\begin{definition}
The active set $\mathcal{A}(x)$ at any feasible $x$ consists of the equality constraints from $\mathcal{E}$ together with the inequality constraints satisfying $c_i (x) = 0$, or equivalently
\[
\mathcal{A}(x) = \mathcal{E} \cup \{i\in\mathcal{I}\ \left|\ c_i(x) = 0\right.\} .
\]
\end{definition}
This is formalized in Lean as bellow.
\begin{lstlisting}
def active_set (point : E) : Finset ℕ := τ ∪ σ.filter fun i : ℕ ↦ p.inequality_constraints i point = (0 : ℝ)
\end{lstlisting}
The term $\sigma$.filter contains only the indexes which satisfy $c_i(x) = 0$.
The active set is made up of those constraints that attains value zero at the current poin. Otherwise, the constraint can be satisfied in a local neighborhood, which makes no difference on the original problem. Hence, the linearized feasible directions are determined by these active constraints as follows.
\begin{definition}
    Given a feasible point $x$ and the active constraint set $\mathcal{A}(x)$, the set of linearized feasible directions $\mathcal{F}(x)$ is
\begin{align}
\mathcal{F}(x) = 
\left\{d \left|
\begin{array}{l}
d^\top \nabla c_i(x) = 0, \quad i \in \mathcal{E}, \\
d^\top \nabla c_i(x) \geqslant 0, \quad i \in \mathcal{A}(x) \cap \mathcal{I}.
\end{array}\right.
\right\}.
\label{linearized feasible}
\end{align}
\end{definition}
We formalize this definition in the namespace \lean{Constrained_OptimizationProblem}.
\begin{lstlisting}
def linearized_feasible_directions (point : E) : Set E := {v | 
    (∀ i ∈ τ, ⟪gradient (p.equality_constraints i) point, v ⟫_ℝ = (0 : ℝ)) ∧ ∀ j ∈ σ ∩ (p.active_set point), ⟪gradient (p.inequality_constraints j) point, v⟫_ℝ ≥ (0 : ℝ)}
\end{lstlisting}
Utilizing the namespace, we can reuse the definition for constrained optimization problem \lean{p} simply by \lean{p.linearized_feasible_directions}. This simplifies the introduction of the definitions for constraint qualifications. 

\subsection{Linear Independent Constraint Qualification}
In this subsection, we introduce the commonly used constraint qualification, LICQ. Recall that we establish the geometric optimality conditions for local minima through the tangent cone in Theorem \ref{thm: geometric optimality}. However, obtaining the tangent cone at a specific point is more challenging compared to deriving linearized feasible directions, which consist of equality and inequality constraints expressed through inner products. Consequently, we aim to connect $T_\Omega(x)$ and $\mathcal{F}(x)$ so as to translate the properties of the tangent cone into linearized feasible directions. Constraint qualifications are conditions under which the linearized feasible set $\mathcal{F}(x)$ coincides with the tangent cone $T_\Omega(x)$. These conditions ensure that $\mathcal{F}(x)$, derived by linearizing the algebraic description of the set at $x$, accurately reflects the essential geometric characteristics of the set in the vicinity of $x$, as encoded by $T_\Omega(x)$.

\begin{definition}
    Given the point $x$ and the active set $\mathcal{A}(x)$, the linear independence constraint qualification (LICQ) holds if the set of active constraint gradients $\{\nabla c_i(x), i \in \mathcal{A}(x)\}$ is linearly independent.
\end{definition}
This is formalized as follows.
\newpage
\begin{lstlisting}
def LICQ (point : E) : Prop := LinearIndependent ℝ
    (fun i : p.active_set point ↦ 
    if i.1 ∈ τ then gradient (p.equality_constraints i.1) point
    else gradient (p.inequality_constraints i.1) point)
\end{lstlisting}

LICQ guarantees the identity $\mathcal{F}(x^\ast) = T_\Omega(x^\ast)$. The proof of this equation is divided into the following two parts.
\begin{theorem}\label{thm: LICQ}
    Let $x^\ast$ be a feasible point. 
    \begin{enumerate}
        \item If the constraints are all differentiable at $x^\ast$, then $T_\Omega(x^\ast) \subset \mathcal{F}(x^\ast)$.
        \item If the constraints are all continuously differentiable at $x^\ast$ and the LICQ conditions are satisfied at $x^\ast$, then $\mathcal{F}(x^\ast) \subset T_\Omega(x^\ast)$. Furthermore, we have $\mathcal{F}(x^\ast) = T_\Omega(x^\ast)$.
    \end{enumerate}
\end{theorem}
The first element of Theorem \ref{thm: LICQ} is formalized in Lean as follows. 
\begin{lstlisting}
theorem linearized_feasible_directions_contain_tagent_cone
    {p : Constrained_OptimizationProblem E τ σ} (xf : x ∈ p.FeasSet)
    (diffable : ∀ i ∈ τ, DifferentiableAt ℝ (equality_constraints p i) x)
    (diffable₂ : ∀ i ∈ σ, DifferentiableAt ℝ (inequality_constraints p i) x) : 
    posTangentConeAt p.FeasSet x ⊆ p.linearized_feasible_directions x
\end{lstlisting}
Notably, the theorem remains valid without the need for a constraint qualification. Regarding the second item in Theorem \ref{thm: LICQ}, although this dependency is not explicitly stated, the proof reveals that the dimension of the domain space is a critical factor. Moreover, the first theorem is established within the general framework of Hilbert spaces. However, in this paper, we only consider the second theorem in $\mathbb{R}^n$, which is defined as \lean{EuclideanSpace ℝ (Fin n)} in Mathlib. There is also small difference between the assumption on the differentiability of the constraint functions. The second element of Theorem \ref{thm: LICQ} requires the continuous differentiability of the constraints. We can see that formalization helps to distinguish which assumptions are crucial for the establishment of the theorem. The formalization of the second element of Theorem \ref{thm: LICQ} is given as below.
\begin{lstlisting}
theorem LICQ_linearized_feasible_directions_sub_posTangentCone
    {p : Constrained_OptimizationProblem (EuclideanSpace ℝ (Fin n)) τ σ}
    (x : EuclideanSpace ℝ (Fin n)) (xf : x ∈ p.FeasSet) 
    (conte : ∀ i ∈ τ, ContDiffAt ℝ (1 : ℕ) (equality_constraints p i) x)
    (conti : ∀ i ∈ σ, ContDiffAt ℝ (1 : ℕ) (inequality_constraints p i) x)
    (LIx : p.LICQ x) (hdomain : p.domain = univ) :
    p.linearized_feasible_directions x ⊆ posTangentConeAt p.FeasSet x
\end{lstlisting}
The term \lean{ContDiffAt ℝ n f x} denotes function $f(x)$ is $n$-order continuous differentiable at $x$. The main difficulties of the proof of the theorems above lie in the second item. We outline the key points of the proof as follows for simplicity.
\begin{enumerate}
\item Define the following notations.
\begin{align*}
    c(x) = \begin{bmatrix}
        c_1(x)\\c_2(x)\\ \vdots\\c_m(x)
    \end{bmatrix},\quad A = \begin{bmatrix}
        \nabla c_1(x^\ast)^\top\\\nabla c_2(x^\ast)^\top\\ \vdots\\\nabla c_m(x^\ast)^\top
    \end{bmatrix},\quad x\in \mathbb{R}^n,\ A \in \mathbb{R}^{m \times n}.
\end{align*}
Since LICQ holds, it holds that the matrix $A$ has full row rank, i.e. $\operatorname{rank}(A) = m$. This is formalized as \lean{LICQ\_Axfullrank}. It is obvious that $m \leqslant n$. Hence, the number of the active constraints is smaller than the dimension of the domain space. We formalize this property as \lean{LICQ\_mlen}. Moreover, $c(x)$ is a function from $\mathbb{R}^n$ to $\mathbb{R}^m$ with Jacobian matrix $A$ at $x^\ast$ in the sense of strict differentiability, which means
\begin{align*}
    c(x) - c(y) - A(x - y) = o(x - y)\ \text{as}\ x, y \to x^\ast.
\end{align*}
This can be deduced from the continuous differentiability of $c_i(x)$ at $x^\ast$ and is formalized in Lean as \lean{LICQ\_strictfderiv\_Ax\_elem}.
\item There exists a matrix $Z$ whose columns are basis for the null space of $A$, i.e.
\begin{align*}
    Z \in\mathbb{R}^{n\times (n-m)},\quad Z \ \text{has full column rank},\quad AZ = 0.
\end{align*}
This is formalized in Lean as \lean{LICQ\_existZ}.
\item Choose an arbitrary $d \in \mathcal{F}(x^\ast)$, and suppose that $\{t_n\}_{n=1}^\infty$ is any sequence of positive scalars such that $\lim\limits_{n\to\infty}t_n = 0$. Define the system of equations $R(z,t) : \mathbb{R}^n \times \mathbb{R} \to \mathbb{R}^n $ as
\begin{align}
    R(z, t) = \begin{bmatrix}
        c(z) - tAd\\Z^\top(z - x^\ast - td)
    \end{bmatrix} = \begin{bmatrix}
        c(z) \\Z^\top(z - x^\ast)
    \end{bmatrix} -t \begin{bmatrix}
       A\\Z^\top
    \end{bmatrix} d= t R_1(z) - tMd = \begin{bmatrix}
        0\\0
    \end{bmatrix}.
    \label{eq: R system}
\end{align}

We can verify that the matrix $M$ is an injection, which is formalized as \lean{LICQ\_injM}. Hence $M$ is invertible. $R_1(z)$ has Jacobian matrix $M$ at $x^\ast$ in the sense of strict differentiability. 
\item According to the implicit function theorem, \eqref{eq: R system} has a unique solution $z_k$ for all values of $t_k$ sufficiently small. This is formalized in Lean as \lean{LICQ\_implicit\_f}. Therefore, we get sequences $z_k \to x^\ast$ and $t_k \to 0$ satisfying $R(z_k,t_k)=0$.
\item We have from \eqref{eq: R system} and the definition of $\mathcal{F}(x^\ast)$ that $z_k$ is indeed feasible. For the constraints in the active set, since $d \in \mathcal{F}(x^\ast)$, the feasibility is derived from \eqref{eq: R system} directly. However, classical textbook omits the discussion on the inactive constraints by assuming all the constraints active. We finish this part by taking Taylor's expansion around the point $x^\ast$ and this is given in formalized theorem \lean{LICQ_inactive_nhds}. Hence, the sequence $\{z_k\}$ coincides with the approaching sequence in the definition of the tangent cone.
\item From the fact that $R(z_k, t_k) = 0$ and $M$ is invertible together,  using Taylor's theorem, it holds that 
\begin{align*}
    0 = M(z_k - x^\ast - td) + o(\lVert z_k - x^\ast\rVert) \implies \frac{z_k - x^\ast}{t_k} = d + o\left(\frac{\lVert z_k - x^\ast\rVert}{t_k}\right).
\end{align*}
Hence, $d \in T_\Omega(x^\ast)$. This is formalized in Lean as \lean{comap1}, \lean{comap2} and \lean{LICQ\_tendsto}.
\end{enumerate}

If LICQ and certain continuity conditions hold, we can substitute $T_\Omega(x^\ast)$ in \eqref{eq: geometric 2} with $\mathcal{F}(x^\ast)$ and get the following theorem:
\begin{theorem}
    Let $x^\ast$ be a local minimizer of the constrained problem. If at $x^\ast$ we have:
\begin{enumerate}
    \item the objective function $f $ is differentiable,
    \item the constraints are all continuous differentiable,
    \item LICQ holds,
\end{enumerate}
then there is no descent direction in the linearized feasible cone, or equivalently
\begin{align}
\mathcal{F}(x^\ast) \cap \{d\ |\ d^\top \nabla f(x^\ast) < 0\} = \varnothing,
\end{align}
which is equivalent to
\begin{align}
\left\{d \left|
\begin{array}{l}
d^\top \nabla f(x^*) < 0, \\
d^\top \nabla c_i(x^*) = 0, \quad i \in \mathcal{E}, \\
d^\top \nabla c_i(x^*) \geqslant 0, \quad i \in \mathcal{A}(x^*) \cap \mathcal{I}
\end{array}\right.
\right\} = \varnothing.
\label{eq: geometric LICQ}
\end{align}
\end{theorem}
The equation \eqref{eq: geometric LICQ} is formalized in Lean as
\newpage
\begin{lstlisting}
theorem local_minima_linearized_feasible_directions_LICQ''
    (loc : EuclideanSpace ℝ (Fin n)) (hl : p.Local_Minima loc)
    (hf : Differentiable ℝ p.objective)
    (conte : ∀ i ∈ τ, ContDiffAt ℝ (1 : ℕ) (equality_constraints p i) loc)
    (conti : ∀ i ∈ σ, ContDiffAt ℝ (1 : ℕ) (inequality_constraints p i) loc)
    (LIx : p.LICQ loc) (hdomain : p.domain = univ):
    ¬ ∃ d : EuclideanSpace ℝ (Fin n),
    (∀ i ∈ τ, ⟪gradient (p.equality_constraints i) loc, d⟫_ℝ = 0) ∧
    (∀ j ∈ σ ∩ p.active_set loc, ⟪gradient (p.inequality_constraints j) loc, d⟫_ℝ ≥ 0) ∧ (⟪gradient p.objective loc, d⟫_ℝ < (0 : ℝ)) 
\end{lstlisting}

\subsection{Linear Constraint Qualification}\label{sec: Linear}
The linear constraint qualification states another assumption for the fact that the tangent cone coincides with the linearized feasible cone. The definition gives as below.
\begin{definition}
    Given the point $x$ and the active set $\mathcal{A}(x)$, we say that the linear constraint qualification holds if all the active constraints $c_i(\cdot), \; \forall i \in \mathcal{A}(x)$ are linear functions. Or equivalently, there exists $a_i \in \mathbb{R}^n$ and $b_i\in \mathbb{R}$ such that 
    \begin{align*}
        c_i(x) = a_i ^\top x + b_i \quad \forall i \in \mathcal{A}(x).
    \end{align*}
\end{definition}
The linear constraint qualification describes a situation when we have full knowledge about the constraint functions, and is also formalized under \lean{Constrained_OptimizationProblem} namespace. Although it might be better to call these linear functions as affine functions,  this constraint qualification is often called linear constraint qualification \cite{nocedal1999numerical}. To describe whether a function is linear or not, we give a new description as \lean{IsLinear}, which is defined on an inner product space. 
\begin{lstlisting}
variable {E : Type _} [NormedAddCommGroup E] [InnerProductSpace ℝ E]

def IsLinear (f : E → ℝ) : Prop := ∃ a, ∃ b, f = fun x ↦ (inner x a : ℝ) + b
def LinearCQ (point : E) : Prop :=
  (∀ i ∈ (p.active_set point ∩ τ), IsLinear (p.equality_constraints i)) ∧
  ∀ i ∈ (p.active_set point ∩ σ), IsLinear (p.inequality_constraints i)
\end{lstlisting}
Under the linear constraint qualification, we can also prove that the linearized feasible directions are the same as the tangent cone at the point. The theorem states as follows.
\begin{theorem}\label{thm: Linear CQ}
    Let $x^\ast$ be a feasible point. If the linear constraint qualificatioin holds at $x^\ast$, it holds that $\mathcal{F}(x^\ast) = T_\Omega(x^\ast)$.
\end{theorem}
The formalization of the theorem is given as below.
\begin{lstlisting}
theorem Linear_linearized_feasible_directions_eq_posTangentCone
    (x : EuclideanSpace ℝ (Fin n)) (xf : x ∈ p.FeasSet)
    (conti : ∀ i ∈ σ \ p.active_set x, ContDiffAt ℝ (1 : ℕ) (inequality_constraints p i) x)
    (Lx : p.LinearCQ x) (hdomain : p.domain = univ):
    p.linearized_feasible_directions x = posTangentConeAt p.FeasSet x 
\end{lstlisting}
The proof of the linear constraint qualification is relatively simpler. We also list the key points of this proof.
\begin{enumerate}
    \item From the first element of Theorem \ref{thm: LICQ}, we obtain the linearized feasible directions contain the tangent cone. The remaining task is to prove that the tangent cone includes the linearized feasible directions. 
    \item Let $v \in \mathcal{F}(x^\ast)$. We need to prove $v \in T_\Omega(x^\ast)$. For the given $x^\ast$, there is a positive scalar $\bar{t}$ such that for any $0 \leq t < \bar{t}$, the inactive constraints remain inactive for $x^\ast + t v$, where $v$ denotes any direction in linearized feasible directions. Or equivalently, it holds
    \begin{align*}
        c_i(x^\ast+t v) > 0, \forall i \in \sigma \backslash \mathcal{A}(x^\ast), \forall \; 0 \leq t < \bar{t}, \; v \in \mathcal{F}(x^\ast).
    \end{align*}
    \item The remaining thing is to construct a sequence approaching $x^\ast$. Define $z_k = x^\ast + (\bar{t}/k) v$ as in Definition \ref{def: tangent cone}. Through straightforward computation with linear constraint qualification, we can examine that $z_k$ is feasible.
\end{enumerate}
It is worthy pointing out that we define a slightly different sequence when formalizing the proof. The definition is given as \lean{let z := fun (k : ℕ) ↦ (t_ / (k + 1)) • v}. This prevents division by zero—a corner case often overlooked in natural language proofs that, while mathematically equivalent, enhances the conciseness of the formalization.
\section{Farkas Lemma}\label{sec: Farkas}
From the constraint qualification, we can describe the geometric optimality condition of the constrained optimization problem through \eqref{eq: geometric LICQ}. To simplify the set on the left side of \eqref{eq: geometric LICQ}, we need to use the Farkas lemma. There are many equivalent descriptions of the Farkas lemma. We give the following version as below. 
\begin{theorem}\label{thm: Farkas}
    Consider a cone $K$:
    \begin{align}
         K = \{By + Cw\ |\ y \geqslant 0\},
         \label{eq: Farkas K}
    \end{align}
    where $B \in \mathbb{R}^{n\times m}, C \in \mathbb{R}^{n\times p}$, $y$ and $w$ are vectors of appropriate dimensions, and $y \geqslant 0$ means that all components of $y$ are nonnegative. Given $g \in \mathbb{R}^n$, then either $g \in K$, or else there is a vector $d \in \mathbb{R}^n$ such that
\begin{align}\label{eq: Farkas theorem}
    g^\top d < 0,\quad B^\top d \geqslant 0,\quad C^\top d = 0.
\end{align}
\end{theorem}

The Farkas lemma in Lean takes the following form:
\begin{lstlisting}
theorem Farkas {τ σ : Finset ℕ} {n : ℕ} {a b : ℕ → EuclideanSpace ℝ (Fin n)} {c : EuclideanSpace ℝ (Fin n)} : (∃ (lam : τ → ℝ), ∃ (mu : σ → ℝ), (∀ i, 0 ≤ mu i) ∧ c = Finset.sum univ (fun i ↦ lam i • a i) + Finset.sum univ (fun i ↦ mu i • b i)) ↔ ¬ (∃ z, (∀ i ∈ τ, inner (a i) z = (0 : ℝ)) ∧ (∀ i ∈ σ, inner (b i) z ≥ (0 : ℝ)) ∧ (inner c z < (0 : ℝ)))
\end{lstlisting}
Although Theorem \ref{thm: Farkas} writes in matrix form and the formalization is given in vector form, they are completely equivalent in mathematics. However, the vector form is easier to use when proving the KKT conditions. The proof of Theorem \ref{thm: Farkas} in \cite{nocedal1999numerical} utilizes detailed analysis on components of the vector and the matrix, which relatively brings difficulties for formalization. Hence, we change another way from the aspect of Carathéodory lemma to formalize the proof the Farkas lemma. 

To complete the proof, we first need to use the property that $K$ is closed. Since $w$ can be split into positive and negative parts, it is easy to see that
\begin{align*}
    K = \left\{
    \left.\begin{bmatrix}
        B & C & -C
    \end{bmatrix}
    \begin{bmatrix}
        y\\w^+\\w^-
    \end{bmatrix}\right|
    \begin{bmatrix}
        y\\w^+\\w^-
    \end{bmatrix} \geqslant 0
    \right\}.
\end{align*}
Therefore, without loss of generality we assume that $K$ has the form
\begin{align*}
    K^\ast = \left\{Ax\left|\right.x\geqslant 0 \right\}.
\end{align*}
The closeness of \eqref{eq: Farkas K} can be deduced with the matrix $A$ of the form $\begin{bmatrix}
        B & C & -C
\end{bmatrix}$.
Consequently, we can view $K^\ast$ as the cone generated by the column vectors of $A$. Denote $\operatorname{Cone}({x_1,\dots,x_n}) = \{x | x = \sum_{i=1}^n a_i x_i, a_i \geq 0\}$. It is easy to see that the cone is a convex set. Let $\widetilde{A} = \{a_1, a_2, \dots, a_m\}$ be the set of the column vectors of $A$, then we have $K^\ast = Cone(\widetilde{A})$. Hence, we transform the proof of the closeness of the set $K^\ast$ to the closeness of a cone.
\begin{lemma}\label{lemma: closed}
    Let $X \subset \mathbb{R}^n$ be a finite set of vectors, then $\operatorname{Cone}(X)$ is closed in $\mathbb{R}^n$.
\end{lemma}
To prove the theorems and lemmas, we give a bunch of new definitions in formalization. The definition of cones in Lean is given as below.
\begin{lstlisting}
def cone {n : ℕ} (s : Finset ℕ) (V : ℕ → (EuclideanSpace ℝ (Fin n))) :
    Set (EuclideanSpace ℝ (Fin n)) := 
    (fun x ↦ Finset.sum s (fun i => x i • V i)) '' {x : ℕ → ℝ | ∀ i : ℕ, 0 ≤ x i}
\end{lstlisting}
The term \lean{f '' s} above denotes the image set of the original set $s$ under the function $f$. Hence,  Lemma \ref{lemma: closed} can be formalized as
\begin{lstlisting}
theorem closed_conic {n : ℕ} (s : Finset ℕ) (V : ℕ → (EuclideanSpace ℝ (Fin n))) : IsClosed (cone s V)
\end{lstlisting}
One elegant proof of this is to use the conic version of Carathéodory lemma:
\begin{lemma}\label{lemma: Caratheodory}
    Let $x \in \operatorname{Cone}(X)$ where $X = \{x_1, x_2, \dots, x_m\}$, $x_i \in \mathbb{R}^n$. Then $x \in \operatorname{Cone}(X')$ for some $X' \subset X$ where the vectors in $X'$ are linearly independent.
\end{lemma}
In Lean this lemma reads
\begin{lstlisting}
theorem conic_Caratheodory {n : ℕ}
    (s : Finset ℕ) (V : ℕ → (EuclideanSpace ℝ (Fin n))) :
    ∀ x ∈ cone s V, ∃ τ : Finset ℕ, (τ ⊆ s) ∧ (x ∈ cone τ V) ∧ (LinearIndependent ℝ (V ∘ coe τ)) ∧
    (∀ σ : Finset ℕ, σ ⊆ s → x ∈ cone σ V → τ.card ≤ σ.card) 
\end{lstlisting}
The proof of this lemma follows the proof of the convex hull version of Carathéodory lemma in Mathlib. We first prove that if $x$ is in the cone generated by the finite set $s$ whose elements are not linear independent, then it is in the cone generated by a strict subset of $s$. This is formalized in \lean{mem_conic_erase}.

Lemma \ref{lemma: Caratheodory} extracts linear independent vectors from the vectors which generate the cone. With this lemma, we can prove the closeness in Leamma \ref{lemma: closed} in the following steps.
\begin{enumerate}
    \item Let $K^\ast = \operatorname{Cone}(\widetilde{A})$, $S = \{\widetilde{B} \subset \widetilde{A}\left.\right| \text{vectors in}\ \widetilde{B}\ \text{are linear independent}\}$, $\widehat{K^\ast} = \bigcup_{\widetilde{B}\in S}\operatorname{Cone}(\widetilde{B})$. Obviously it holds that $\widehat{K^\ast} \subset K^\ast$.
    \item From Lemma \ref{lemma: Caratheodory}, $K^\ast$ is generated by the conic combination of a set of linear independent vectors. It yields that $K^\ast \subset \widehat{K^\ast}$, hence $K^\ast = \widehat{K^\ast}$.
    \item Let $B \in \mathbb{R}^{n\times p}$ be the matrix whose column vectors are $\widetilde{B}$. Since vectors in $\widetilde{B}$ are linear independent, the linear map $f(x) = Bx$ is an injection. Therefore, $f(x)$ is a closed embedding from $\mathbb{R}^p$ to $\mathbb{R}^n$. Note that $\operatorname{Cone}(\widetilde{B}) = f(\Omega)$, where $\Omega = \{x\in \mathbb{R}^p\left.\right| x\geqslant 0\}$ is closed in $\mathbb{R}^p$. Since a closed embedding is also a closed map, we obtain the closeness of $f(\Omega)$. This is formalized in Lean as \lean{closed\_conic\_idp}.
    \item $\widehat{K^\ast}$ is a finite union of closed sets and is therefore also closed. This is formalized in Lean as \lean{cone\_eq\_finite\_union} and \lean{closed\_conic}.
\end{enumerate}

The remaining part of Theorem \ref{thm: Farkas} involves discussing the two cases. The most important is to use the hyperplane separation theorem of the convex set. Prepared with the fact that the set $K$ in \eqref{eq: Farkas K} is a closed and convex set, we prove the statements as follows.
\begin{enumerate}
    \item If we have $g \in K$, it holds that there exist $y \geq 0$ and $w \in \mathbb{R}^p$ satisfying $g = By + Cw$. Hence it yields 
    \begin{align*}
        g^\top d = y^ \top B^\top d + w^\top C^\top d \geq 0,
    \end{align*}
    which is in contradiction to $g^\top d <0$.
    \item If there is no vector $d$ satisfying equation \eqref{eq: Farkas theorem}, we prove the target result by contradiction. Hence $g \notin K$. From the hyperplane separation theorem of the convex set, which can be viewed as a geometric form of the Hahn-Banach theorem, we can get there exists a vector $d$ and a constant $\alpha$ satisfying 
    \begin{align} \label{eq: separation}
        d^ \top g < \alpha < d^\top z, \quad \forall z \in K.
    \end{align}
    The theorem needed is well formalized in Mathlib as \lean{geometric_hahn_banach_point_closed}. The fact that $K$ is closed, which we pay much effort to prove, is used as an assumption for this theorem. Notice that $0 \in K$, therefore $g^\top d < 0$. Moreover, for any column vector of matrix $B$ as $b_i$, it holds that $tb_i \in K$, where $t$ is a nonnegative scalar. Hence, we have $t d^\top b_i > \alpha$ for any $t \geq 0$. Letting $t$ go to infinity, it yields $d^\top b_i \geq 0$. Consequently, $B^\top d \geq 0$. We can apply similar analysis on the matrix $C$ and get $C^\top d = 0$. A vector $d$ satisfying \eqref{eq: Farkas theorem} is constructed, which gives a contradiction.
\end{enumerate}

The discussion is fully unfolded in formalization. Proves omitted for writing simplicity are replenished to ensure there is no mistake. More general forms of Farkas lemma are also considered in the literature, such as the conic Farkas lemma. One of the key points of the proof is to check whether the cone $K$ is a closed set or not. This is not true when considering general cones rather than $\mathbb{R}^n_+$. The formalization of these situations is beneficial to the formalization of general KKT conditions, which we put to future work.   

\section{KKT Conditions}\label{sec: KKT}

The KKT conditions of constrained optimization problem \eqref{eq: constrained problem} provides a necessary criterion for a local minimizer. Based on the Lagrangian function \eqref{eq: lagrangian function}, the KKT conditions state the following.
\begin{theorem}
\label{thm: KKT}
    Let $x^\ast$ be a local minimizer of the constrained problem that $f $ is differentiable at $x^\ast$, all constraints are continuous differentiable at $x^\ast$, and either linear independent constraint qualification or linear constraint qualification holds at $x^\ast$. Then there are Lagrange multiplier vectors $\lambda^\ast$ and $\mu^\ast$, such that the following conditions are satisfied at $(x^\ast, \lambda^\ast, \mu^\ast)$:
\begin{align*}
    \nabla_x \mathcal{L}(x^\ast, \lambda^\ast, \mu^\ast) &= 0, & \text{(Stationarity)}\\
    c_i(x^\ast) &= 0,\quad \forall i \in \mathcal{E}, & \text{(Primal Feasibility)}\\
    c_i(x^\ast) &\geqslant 0,\quad \forall i \in \mathcal{I}, & \text{(Primal Feasibility)}\\
    \mu_i^\ast &\geqslant 0, \quad \forall i \in \mathcal{I}, & \text{(Dual Feasibility)}\\
    \mu_i^\ast c_i(x^\ast) &= 0, \quad \forall i \in \mathcal{I}. & \text{(Complementary Slackness)}
\end{align*}
\end{theorem}
The Lagrange multipliers corresponding to inactive inequality constraints are zero through complementary slackness. The proof of the KKT condition in Theorem \ref{thm: KKT} is straightforward based on the preparation in this paper. Using the Farkas lemma, we deduce that equation \eqref{eq: geometric LICQ} is equivalent to
\begin{align}
    \nabla f(x^\ast) = \sum_{i\in\mathcal{A}(x^\ast)}\mu_i^\ast\nabla c_i(x^\ast),\quad \mu_i^\ast\geqslant 0 \ \text{for any }\ i \in \mathcal{A}(x^\ast)\cap\mathcal{I}.
    \label{eq: kkt grad1}
\end{align}
Define the vector $\mu^\ast$ to be
\begin{align*}
    \mu_i^\ast = \begin{cases}
        \mu_i,\quad&i \in \mathcal{A}(x^\ast),\\
        0,\quad & i \in \mathcal{I}\backslash\mathcal{A}(x^\ast).
    \end{cases}
\end{align*}
Then equation \eqref{eq: kkt grad1} can be written as
\begin{align}
    \nabla f(x^\ast) = \sum_{i\in\mathcal{I}\cup\mathcal{E}}\mu_i^\ast\nabla c_i(x^\ast),\quad \mu_i^\ast\geqslant 0\ \text{and}\ \mu_i^\ast c_i(x^\ast) = 0 \ \text{for}\ i \in \mathcal{I}.
    \label{eq: kkt grad2}
\end{align}
Thus we get the stationarity, dual feasibility and complementary slackness conditions in KKT. At this point, we have obtained all the KKT conditions in \eqref{thm: KKT}. In Lean, Theorem \ref{thm: KKT} with linear independent constraint qualification is formalized as follows.
\begin{lstlisting}
theorem first_order_neccessary_LICQ
    (p1 : Constrained_OptimizationProblem (EuclideanSpace ℝ (Fin n)) τ σ)
    (loc : EuclideanSpace ℝ (Fin n)) (hl : p1.Local_Minima loc)
    (hf : Differentiable ℝ p1.objective)
    (conte : ∀ i ∈ τ, ContDiffAt ℝ (1 : ℕ) (equality_constraints p1 i) loc)
    (conti : ∀ i ∈ σ, ContDiffAt ℝ (1 : ℕ) (inequality_constraints p1 i) loc)
    (hLICQ : p1.LICQ loc) (hdomain : p1.domain = univ) :
    p1.FeasPoint loc ∧
    (∃ (lambda1 : τ → ℝ) (lambda2 : σ → ℝ), (gradient (fun m ↦ (p1.Lagrange_function m lambda1 lambda2)) loc = 0) ∧
    (∀ j : σ, lambda2 j ≥ 0) ∧
    (∀ j : σ, (lambda2 j) * (p1.inequality_constraints j loc) = 0))
\end{lstlisting}
Similar formalization results can be deduced for linear constraint qualification. The difference is that we replace the assumption \lean{p1.LICQ loc} with \lean{p1.LinearCQ loc}.
\begin{lstlisting}
theorem first_order_neccessary_LinearCQ
    (p1 : Constrained_OptimizationProblem (EuclideanSpace ℝ (Fin n)) τ σ)
    (loc : EuclideanSpace ℝ (Fin n)) (hl : p1.Local_Minima loc)
    (hf : Differentiable ℝ p1.objective)
    (conti : ∀ i ∈ σ \ p1.active_set loc, ContDiffAt ℝ (1 : ℕ) (inequality_constraints p1 i) loc)
    (hLinearCQ : p1.LinearCQ loc) (hdomain : p1.domain = univ) :
    p1.FeasPoint loc ∧ (∃ (lambda1 : τ → ℝ) (lambda2 : σ → ℝ),
    (gradient (fun m ↦ (p1.Lagrange_function m lambda1 lambda2)) loc = 0) ∧ (∀ j : σ, lambda2 j ≥ 0) ∧
    (∀ j : σ, (lambda2 j) * (p1.inequality_constraints j loc) = 0)) 
\end{lstlisting}
Apart from these properties, we can formalize the general KKT point for a constrained optimization problem as follows.
\begin{lstlisting}
def KKT_point (self : Constrained_OptimizationProblem E τ σ)
    (x : E) (lambda : τ → ℝ) (mu : σ → ℝ) : Prop :=
  (gradient (fun m ↦ (self.Lagrange_function m lambda mu)) x = 0) 
  ∧ (x ∈ self.FeasSet) ∧ (mu ≥ 0) 
  ∧ (∀ i : σ, mu i * self.inequality_constraints i x = 0)
\end{lstlisting}

The formalization of KKT conditions presented is critical not only for optimization theory, but also important for identifying local minima in concrete optimization problems. Based on the formalization, we can deduce the optimality conditions for an optimization problem with certain properties on the constraints and objective function.

\section{Duality Property}\label{sec: weak duality}
In this subsection, we introduce the duality properties of the primal problem \eqref{eq: constrained problem}. From the dual problem, we can easily get a lower bound for the primal problem, which is derived from the weak duality property and will be formalized in this subsection. In some cases, the dual problem obtains the same optimal value as the primal problem. Since the dual problem may be easier to solve than the primal problem, some algorithms seek to solve the dual problem instead of the primal problem. The definition of the dual problem starts by introducing the dual objective function $q(\lambda, \mu)$ from the Lagrangian function of the primal problem \eqref{eq: lagrangian function} as:
\begin{align*}
    q(\lambda, \mu) \triangleq \inf\limits_x \mathcal{L}(x, \lambda, \mu).
\end{align*}
The definition is also formalized under the namespace \lean{Constrained_OptimizationProblem}. We use the same notations as in the formalization the KKT conditions. 
\begin{lstlisting}
variable {E : Type _} {τ σ : Finset ℕ}
variable {p : Constrained_OptimizationProblem E τ σ}
def dual_objective : (τ → ℝ) → (σ → ℝ) → EReal :=
  fun (lambda1 : τ → ℝ) (lambda2 : σ → ℝ) ↦
    ⨅ x ∈ p.domain, (p.Lagrange_function x lambda1 lambda2).toEReal
\end{lstlisting}
For some $(\lambda, \mu)$ pair, the infimum might become infinity. In descriptions using natural language words, we often overlook these pairs and take $q(\lambda, \mu)$ as a real-valued function \cite{nocedal1999numerical}. However, we cannot do this in formalization. The negative infinity is taken as $\bot$ in the type \lean{EReal}, but not defined in the type $\mathbb{R}$. The functions \lean{Real.toEReal} and \lean{EReal.toReal} help convert between these two types. The former takes the natural embedding between real numbers and extended real numbers. When transforming the bottom and top elements of \lean{EReal} to $\mathbb{R}$, they are mapped to $0$. Hence, in the definition of the dual objective function, we convert the value of the Lagrangian function to \lean{EReal}, and then take the infimum over the domain set of the primal variable. 

Based on the dual objective function, we can give the definition of the dual problem toward the primal problem. Note that the objective function in the constrained optimization problem \eqref{eq: constrained problem} is a real-valued function. Therefore, we need another conversion. The dual problem gives as:
\begin{align}\label{eq: dual problem}
    \max\limits_{\lambda, \mu} q(\lambda, \mu)  \quad  \text{subject to}\  \quad \mu \geq 0.
\end{align}
The domain of this problem is all the $(\lambda, \mu)$ such that the dual objective function does not attain negative infinity. The formalization gives as:
\begin{lstlisting}
def dual_problem : Constrained_OptimizationProblem ((τ → ℝ) × (σ → ℝ)) ∅ σ where
  domain := {x | p.dual_objective x.1 x.2 ≠ ⊥}
  equality_constraints := fun _ _ ↦ 0
  inequality_constraints := fun i x ↦ if h : i ∈ σ then x.2 ⟨i, h⟩ else 0
  objective := fun x ↦ (p.dual_objective x.1 x.2).toReal
  eq_ine_not_intersect  simp
\end{lstlisting}
In \eqref{eq: dual problem}, there is no equality constraint. Hence, in the formal definition, the indices of the equality constraint are the empty set. It is natural to ask whether the dual objective function can attain positive infinity or not. Assume that the domain of the primal problem is not empty, which is a trivial assumption on the primal problem. Then we have $(\lambda, \mu)$, $q(\lambda, \mu) < \infty$. Note that we only require the domain set of the primal problem to be nonempty. The feasible set may be empty due to the constraints of the primal problem. This property is formalized as follows.
\begin{lstlisting}
lemma dual_objective_le_top {p : Constrained_OptimizationProblem E τ σ} (hp : (p.domain).Nonempty) :
    ∀ lambda1 lambda2, p.dual_objective lambda1 lambda2 < ⊤ 
\end{lstlisting}

We are ready to present the weak duality property. The weak duality property bridges the supremum value of the dual problem between the infimum value of the primal problem.
\begin{theorem}\label{thm: weak duality}
    For any $x$ feasible for \eqref{eq: constrained problem} and any $\lambda$, $\mu \geq 0$, it holds that $q(\lambda, \mu) \leq f(x)$. Or equivalently, the infimum value of the primal problem is larger than the supremum value of the dual problem.
\end{theorem}
We first formalize the infimum value and the supremum value of constrained optimization problems. These two definitions are also in the namespace \lean{Constrained_OptimizationProblem}.
\begin{lstlisting}
def inf_value : EReal := sInf (Real.toEReal '' (p.objective '' (p.FeasSet)))
def sup_value : EReal := sSup (Real.toEReal '' (p.objective '' (p.FeasSet)))
\end{lstlisting}
The terms \lean{sInf} and \lean{sSup} denote the infimum and supremum of a set respectively in Lean. Another conversion between $\mathbb{R}$ and \lean{EReal} is utilized to ensure the value can reach infinity. Slightly different from the natural language in Theorem \ref{thm: weak duality}, we avoid using $f(x)$ to consider the case where the feasible set of the primal problem is empty. Though this is relatively trival, it contains some corner case. We provide two formalization version of the theorem
\begin{lstlisting}
theorem weak_duality {p : Constrained_OptimizationProblem E τ σ}
    (lambda1 : τ → ℝ) {lambda2 : σ → ℝ} (mpos : lambda2 ≥ 0):
    p.dual_objective lambda1 lambda2 ≤ p.inf_value 
\end{lstlisting}

\begin{lstlisting}
theorem weak_duality' {p : Constrained_OptimizationProblem E τ σ} :
    (p.dual_problem).sup_value ≤ p.inf_value 
\end{lstlisting}
The proof the weak duality property mainly involves the following inequality. For any $\lambda$ and $\mu \geq 0$ and any feasible $\bar{x}$, we have 
\begin{align}\label{eq: weak duality proof}
    q(\lambda, \mu) &= \inf\limits_x f(x) - \sum_{i \in \mathcal{E}} \lambda_i c_i(x) - \sum_{i \in \mathcal{I}} \mu_i c_i(x) 
    \leq  f(\bar{x}) - \sum_{i \in \mathcal{E}} \lambda_i c_i(\bar{x}) - \sum_{i \in \mathcal{I}} \mu_i c_i(\bar{x}) \leq f(\bar{x}).
\end{align}
The final inequality follows from the fact that $\bar{x}$ is feasible and $\mu \geq 0$. Taking infimum over $\bar{x}$ in \eqref{eq: weak duality proof}, we can obtain Theorem \ref{thm: weak duality}.

Based on the weak duality, it is natural to ask when the solution of the primal problem and the dual problem has a connection. This is studied under the assumption of convexity. We consider constrained optimization problems without equality constraints for simplicity. The form of the problem is given as below
\begin{align}\label{eq: weak duality problem}
    \min\limits_x f(x) , \quad \text{subject to}\  \quad c(x) \geq 0,
\end{align}
where $c(x) = (c_1(x), \cdots, c_m(x))^\top \in \mathbb{R}^m$ denotes the inequality constraints. The Lagrangian function of problem \eqref{eq: weak duality problem} is $\mathcal{L}(x,\mu) = f(x) - \mu^\top c(x)$. The KKT conditions of problem \eqref{eq: weak duality problem} are given as below:
\begin{align}\label{eq: KKT weak duality}
\begin{aligned}
    \nabla_x \mathcal{L}(x,\mu) = \nabla f(x) - \nabla c(x) \mu &= 0, \\
    c(x) &\geq 0,\\
    \mu &\geq 0, \\
    \mu_i c_i(x) &= 0, \quad i=1,\cdots, m.
\end{aligned}
\end{align}
We make the following basic assumptions of  problem \eqref{eq: weak duality problem}. 
\begin{assumption} \label{assump: weak duality}
    \begin{enumerate}
        \item The objective function $f(x)$ is convex, and all the inequality constraints $c_i(x)$ are concave.
        \item The objective function $f(x)$ and all the inequality constraints $c_i(x)$ are differentiable at the solution of the primal problem.
    \end{enumerate}
\end{assumption}
Under this assumption, we can establish that optimal Lagrange multipliers for \eqref{eq: weak duality problem} are solutions of the dual problem under certain conditions.
\begin{theorem}\label{thm: KKT dual}
    Suppose that $\bar{x}$ is an optimal solution of problem \eqref{eq: weak duality problem} and that Assumption \ref{assump: weak duality} holds. Then for any $\bar{\mu}$ for which $(\bar{x}, \bar{\mu})$ satisfies \eqref{eq: KKT weak duality} is a solution of the problem \eqref{eq: weak duality problem}.
\end{theorem}
This theorem is formalized as follows.
\begin{lstlisting}
theorem KKT_multipliers_solution_dual_problem {p : Constrained_OptimizationProblem E τ σ}
    (h : ConvexOn ℝ univ p.objective) (hτ : τ = ∅) (x : E)
    (hf : DifferentiableAt ℝ p.objective x) 
    (hi : ∀ i ∈ σ, ConcaveOn ℝ univ (inequality_constraints p i))
    (conti : ∀ i ∈ σ, DifferentiableAt ℝ (inequality_constraints p i) x)
    (lambda1 : τ → ℝ) (lambda2 : σ → ℝ) (hg : p.Global_Minimum x)
    (hKKT : KKT_point p x lambda1 lambda2) :
    (p.dual_problem).Global_Maximum (lambda1, lambda2) 
\end{lstlisting}
In the formalization, we use the assumption \lean{hτ} to show that there is no equality constraint. The proof consists of the following key points.
\begin{enumerate}
    \item From the first element in Assumption \ref{assump: weak duality}, the Lagrangian function is convex, which is formalized in \lean{convex_problem_convex_Lagrange}. From the first order optimality condition of the convex functions, it yields for any $x$
    \begin{align*}
        \mathcal{L}(x, \bar{\mu}) \geq \mathcal{L}(\bar{x}, \bar{\mu}) + \nabla_x \mathcal{L}(\bar{x}, \bar{\mu}) ^\top (x-\bar{x}) = \mathcal{L}(\bar{x}, \bar{\mu}).
    \end{align*}
    Hence we can get that the infimum of $\mathcal{L}(x, \bar{\mu})$ is attained at $(\bar{x}, \bar{\mu})$.
    \item Therefore, it holds that
    \begin{align*}
        q(\bar{\mu}) = \inf_x \mathcal{L}(x, \bar{\mu}) = \mathcal{L}(\bar{x}, \bar{\mu}) = f(\bar{x}) -\mu^\top c(\bar{x}) = f(\bar{x}).
    \end{align*}
    The last equality comes from the complementary slackness of the KKT conditions. From the weak duality \eqref{eq: KKT weak duality}, we have for any $q(\mu) \leq f(\bar{x})$. Hence $\bar{\mu}$ is a solution of the dual problem \eqref{eq: weak duality problem}. 
\end{enumerate}
\section{Conclusion and Future Work}\label{sec: conclution}
In this paper, we discuss the formalization of optimality conditions for smooth constrained optimization problems. First, we give the formalization of the basic definitions in constrained optimization problems, including the objective functions, the constraint functions, the local minima and the Lagrangian function. To prove the KKT conditions, we start from the geometric optimality conditions. Then we formalize the constraint qualifications. Both linear independent constraint qualification and the linear constraint qualification are considered to prove the equivalence between the tangent cone and the linearized feasible direction. The Farkas lemma is used to further analyze the linearized feasible directions. Based on these results, we give the formalization of the KKT conditions. The dual problem and duality property are also formalized. 

The study of the optimality conditions has much more to be formalized and these can be considered in future work. For constraint qualifications, many other weaker or equivalent constraint qualifications are considered to deal with different situations in literature, such as the Mangasarian–Fromovitz constraint qualification (MFCQ). For optimality conditions, we restrict the discussion within the first order optimality conditions in this paper. However, the KKT conditions are only necessary conditions for local minima. Second-order optimality conditions provide necessary and sufficient criteria for local minima. In the context of convex optimization problems, strong duality is guaranteed by the Slater condition, which ensures the existence of interior points in the feasible set.




\bibliography{KKT-formalization}

\end{document}